\documentclass[journal,twoside,web]{ieeecolor}

\usepackage{generic}
\usepackage{cite}
\usepackage{amsmath,amssymb,amsfonts,color}
\usepackage{algorithmic}
\usepackage{graphicx}
\usepackage{textcomp}
\usepackage{orcidlink}
\usepackage{subfigure}
\usepackage{epsfig} 
\usepackage{cite}
\usepackage{lcsys}
\usepackage{hyperref}
\hypersetup{colorlinks,allcolors=black}

\def\BibTeX{{\rm B\kern-.05em{\sc i\kern-.025em b}\kern-.08em
    T\kern-.1667em\lower.7ex\hbox{E}\kern-.125emX}}
\newtheorem{theorem}{Theorem}[section]

\newtheorem{definition}[theorem]{Definition}
\newtheorem{remark}[theorem]{Remark}

\newcommand{\lr}[2]{\langle #1, #2\rangle}
\begin{document}

\title{Computational Optimal Transport and Filtering on Riemannian Manifolds\thanks{Approved for public release; distribution is unlimited. Public Affairs approval $\#$ AFRL-2023-4936. The views expressed are those of the authors and do not necessarily reflect official policy or position of the Department of the Air Force, the Department of Defense or the U.S. Government. }}

\author{Daniel Grange$^{*}$\orcidlink{0000-0002-7643-496X} \and 
Mohammad Al-Jarrah$^{\dagger}$\orcidlink{0009-0006-0433-9230}  \and 
\and Ricardo Baptista$^{\ddagger}$\orcidlink{0000-0002-0421-890X} \and Amirhossein Taghvaei$^{\dagger}$\orcidlink{0000-0002-1536-892X}  \and\\ 
Tryphon T.\ Georgiou$^{\star}$\orcidlink{0000-0003-0012-5447} \and
Sean Phillips$^{\S}$\orcidlink{/0000-0001-9074-6049} \and
Allen Tannenbaum$^{*}$\orcidlink{0000-0002-0567-5256}
\thanks{$^{*}$Department of Computer Science, Stony Brook University, NY; \texttt{\{daniel.grange,allen.tannenbaum\}@stonybrook.edu}.}
\thanks{$^{\dagger}$Department of Aeronautics and Astronautics, University of Washington, Seattle, WA; \texttt{\{mohd9485,amirtag\}@uw.edu}.}
\thanks{$^{\ddagger}$Department of Computing and Mathematical Sciences, California Insitute of Technology, Pasadena, CA; \texttt{rsb@caltech.edu}.}
\thanks{$^{\star}$Department of Mechanical and Aerospace Engineering, University of California, Irvine, CA; \texttt{tryphon@uci.edu}.}
\thanks{$^{\S}$ Space Vehicles Directorate, Air Force Research Laboratory, Kirtland AFB, NM}
\thanks{DG and MAJ led the computational effort of the project, A.\ Taghvaei led the methodological effort, largely built upon \cite{al2023optimal}, and the remaining authors supervised the development of the work. The research was supported in part by the  NSF under grant EPCN-2318977, AFOSR under FA9550-23-1-0096 and ARO under W911NF-22-1-0292.}
}
\maketitle
\thispagestyle{empty}

\begin{abstract}
   In this paper we extend recent developments in computational optimal transport to the setting of Riemannian manifolds. In particular, we show how to learn optimal transport maps from samples that relate probability distributions defined on manifolds. Specializing these maps for sampling conditional probability distributions provides an ensemble approach for solving nonlinear filtering problems defined on such geometries. The proposed computational methodology is illustrated with examples of transport and nonlinear filtering on Lie groups, including the circle $S^1$, the special Euclidean group $SE(2)$, 
   and the special orthogonal group $SO(3)$.
\end{abstract}

\begin{IEEEkeywords}
Optimal Transportation, Optimal Control, Nonlinear Filtering, Riemannian manifolds. 
\end{IEEEkeywords}

\section{Introduction}
\label{sec:introduction}
\IEEEPARstart{T}{he} theory of optimal transport (OT) has emerged as a powerful mathematical tool in a wide range of engineering and control applications \cite{sepulchre2021optimal}. This is largely due to the fact that it induces a natural and computationally tractable geometry on the space of probability distributions~\cite{Vil08,peyre2019computational}. 
The metric that the theory provides  to quantify distance between distributions, the {\it Wasserstein metric}, gives rise to natural geodesic flows and transport maps that can be used to interpolate, average, and correspond distributions in a physically meaningful sense. For these reasons, OT has proven enabling for an ever expanding range of applications in machine learning~\cite{arjovsky2017wasserstein,courty2015optimal,peyre2019computational}, and image processing~\cite{kolouri2017optimal,dominitz2010texture,rabin2011wasserstein,su2015optimal}, besides ones in control and estimation~\cite{sepulchre2021optimal,chen2021optimal,chen2021optimal2,taghvaei2021OTFPF}.  

In this rapidly developing landscape of OT techniques and applications, neural networks and stochastic optimization have come to provide a potentially transformative framework for the development of efficient and scalable numerical algorithms~\cite{CheGeoTan18,leygonie2019adversarial,xie2019scalable,makkuva2020optimal,korotin2021wasserstein}. The focus so far on utilizing such techniques however has been limited to applications of OT on Euclidean spaces.
Yet, optimal transport can equally well be considered on manifolds, that are especially relevant in control and robotic applications. A manifold structure is naturally imposed by geometric constraints, as in attitude estimation of  aircraft~\cite{hua2013implementation,barrau2014intrinsic}, localization of mobile robots~\cite{barczyk2015invariant,hesch2014camera}, 
and visual tracking of humans and objects~\cite{kwon2013geometric,choi2012robust}.
%

Thus, one of the goals of the present paper is to 
develop a computational framework for OT in the setting of Riemannian manifolds~\cite{mccann2001polar,ambrosio2005gradient} with special attention to matrix Lie-groups, as these encompass the majority of the motivating applications. A second goal of the paper is to use the framework for sampling conditional distributions, in order to perform nonlinear filtering  on Riemannian manifolds. 

%
%
Specifically, 
we make the following key contributions:
\begin{itemize}
    \item[i)] We propose a sample-based 
    computational methodology for computing OT maps on Riemannian manifolds. Our proposed methodology extends the min-max formulation of~\cite{makkuva2020optimal} by combining it with McCann's characterization of optimal transport maps~\cite{mccann2001polar}.
    \item[ii)] We propose a sample-based and likelihood-free method to sample conditional distributions on Riemannian manifolds. In order to do so, we use the recently introduced framework of block-triangular transport maps, that is used in the context of conditional generative models~\cite{kovachki2020conditional,shi2022conditional} and nonlinear filtering~\cite{spantini2019coupling,taghvaei2022optimal}. 
    \item[iii)] We illustrate our proposed algorithms on several numerical examples on the circle, special Euclidean group $SE(2)$, and the special orthogonal group $SO(3)$.
\end{itemize}





\section{Problem formulation and background}
Let $\mathcal M$ be a smooth connected manifold without boundary that is equipped with the Riemannian metric $\lr{\cdot}{\cdot}_g$.  Let $d(z,z')$ denote the geodesic distance for any $z,z'\in M$. 
We are interested in solving the following two  problems. 

\noindent
{\bf Optimal control:}
This is the problem to steer a random process $X(t)$, taking values in $\mathcal M$, from an initial probability distribution $P$ to a terminal probability distribution $Q$. It is formulated as follows: 
\begin{equation}\label{eq:opt-control}
\begin{aligned}
    &\min_u\, \mathbb E\left[\int_0^1 \|u(t)\|_g^2\,dt\right],\\
    &\text{s.t}\quad \dot X(t) = u(t),\quad X(0) \sim P,\quad X(1)\sim Q,
\end{aligned}
\end{equation}
where the control input $u(t) \in T_{X(t)}\mathcal M$ for all $t\in[0,1]$, and the control cost is the square of the Riemannian norm $\|u\|_g := \sqrt{\lr{u}{u}_g}$. In practice, such problems arise when controlling an ensemble of agents or a swarm of robots.   




\noindent
{\bf Optimal filtering:}
The second problem we are interested is to compute the conditional distribution of a hidden random variable $X \in \mathcal M$ given  
 an observed random variable $Y \in \mathbb R^m$.  The conditional distribution of $X|Y$, i.e., the posterior distribution, is given by Bayes' law as
\begin{equation}\label{eq:Bayes}
     P_{X|Y}(x|y) = \frac{P_{Y|X}(y|x)P_X(x)}{P_Y(y)},
 \end{equation}
 where $P_X$ is the prior probability distribution of $X$, $P_{Y|X}$ is the likelihood of observing $Y$ given $X$, and  $P_Y(y) = \int_{\mathcal M} P_{Y|X}(y|x)P_X(x) dx$ is the probability distribution of $Y$. Sampling the conditional distribution in~\eqref{eq:Bayes} is an essential step in many nonlinear ensemble filtering algorithms~\cite{doucet09}. Classic algorithms include particle filters (or sequential Monte Carlo methods), which suffer from weight degeneracy ~\cite{cheng2010particle, vernaza2006rao}, and Kalman-filter-type algorithms, which often fail to represent multi-modal distributions~\cite{crassidis2003unscented,bonnabel2009IEKF,barrau2015TAC, barczyk2015invariant}.   
 
In Section~\ref{sec:methodology} we will show that, in general, both problems can be formulated as problems of OT on Riemannian manifolds. But before we proceed, we review next some key results of the theory of OT on Riemannian manifolds. 
 
\subsection{Background on OT on Riemannian manifolds} 
Given two 
probability distributions $P$ and $Q$ on $\mathcal M$, the Monge optimal transportation problem seeks a map $T\colon \mathcal M\to \mathcal M$ that solves the optimization problem
\begin{align}\label{eq:Monge}
	\inf_{T\in \mathcal T(P,Q)}\, \mathbb E_{Z\sim P}[c(Z,T(Z))],
\end{align}
where $\mathcal T(P,Q):=\{T\colon \mathcal M\to \mathcal M;\,T_{\#}P=Q\}$ is the set of all transport maps pushing forward $P$ to $Q$, and $c\colon \mathcal M \times \mathcal M \to \mathbb R$ is a lower semi-continuous cost function that is bounded from below.  To account for the challenging nonlinear constraint in~\eqref{eq:Monge}, the Monge problem is relaxed by replacing deterministic transport maps with stochastic couplings and solving
\begin{align}\label{eq:Kantorovich}
	\inf_{\pi \in \Pi(P,Q)}\, \mathbb E_{(Z,Z')\sim \pi}[c(Z,Z')],
\end{align}
instead, where $\Pi(P,Q)$ denotes the set of all joint distributions on $\mathcal M \times \mathcal M$ with marginals  $P$ and $Q$. This relaxation, due to Kantorovich, turns the Monge problem into a linear program, whose dual becomes
\begin{align}\label{eq:dual-Kantorovich}
	\sup_{(\phi,\psi)\in \text{Lip}_c}\,\mathbb E_{Z\sim P}[\phi(Z)] + E_{Z'\sim Q}[\psi(Z')],
\end{align}
where $\text{Lip}_c$ is the set of pairs  of functions $\phi,\psi$ from $\mathcal M\to\mathbb R$, 
\begin{align*}
\text{Lip}_c&:=\{(\phi,\psi)\mid\,\phi(z)+\psi(z')\leq c(z,z'),\,\forall z,z'\in \mathcal M\}.
\end{align*}

\begin{definition}
	Given a function $\phi:\mathcal M \to \mathbb R \cup \{\pm \infty\}$, the inf-$c$ convolution is given by 
	\begin{align*}
		\phi^c(z') = \inf_{z\in M}\left[c(z,z')-\phi(z)\right]. 
	\end{align*}
	Moreover, $\phi$ is said to be $c$-concave if $\exists\, \psi$ such that $\phi=\psi^c$. 
\end{definition}


\begin{theorem}[McCann~\cite{mccann2001polar}] Consider the Monge problem~\eqref{eq:Monge} with cost $c(z,z')=d(z,z')^2/2$, and assume that $P$ is absolutely continuous with respect to the volume measure on $\mathcal M$. Then, there exists a unique minimizer $T(z) = \exp_z[-\nabla \phi(z)]$ where $\phi$ is $c$-concave and the pair $(\phi,\phi^c)$ maximize the dual Kantorovich problem~\eqref{eq:dual-Kantorovich}.
 \label{thm:McKann}
\end{theorem}
\medskip 

\begin{remark}
	For the case where $\mathcal M=\mathbb R^n$ and $c(z,z')=\frac{1}{2}\|z-z'\|_2^2$, the optimal map is given by $T(z)=z-\nabla  \phi(z)$, where $\|z\|^2_2/2-\phi(z)$ is a convex function; this is a celebrated result due to Brenier~\cite{brenier1991polar}.
\end{remark}

\subsection{Computational methods for OT on manifolds}
The majority of existing computational algorithms for OT on manifolds are concerned with constructing normalizing flows and $c$-concave potential functions~\cite{cohen2021riemannian,rezende2020normalizing,lou2020neural,mathieu2020riemannian}. For instance, \cite{cohen2021riemannian} proposes a general method to represent $c$-concave functions by taking the inf-$c$ convolution of a piecewise constant function, reference \cite{rezende2020normalizing} introduces a set of diffeomorphisms on circles, tori, and spheres that are  used as building blocks of a neural network to represent $c$-concave functions, and 
reference \cite{de2022riemannian} is concerned with implementation of diffusion models on Riemannian manifolds via a variety of approaches that are based on projection, using a Lie-algebra basis and coordinate vector-fields, in order to represent general vector-fields. Lastly, \cite{mathieu2020riemannian} introduces the geodesic distance layer, which generalizes the concept of linear layers to
manifolds, used as an input layer to a neural net to represent vector-fields on manifolds. 

In the present paper, we also use neural networks to represent transport maps on 
the given manifold. Our neural network architectures, along with the special details for each example that we present, are explained in Section~\ref{sec:numerics}.
 \section{Solution methodology} \label{sec:methodology}


In this section, we present the OT formulation of the two problems, optimal control and filtering, along with a stochastic optimization formulation for their numerical solution. 
\subsection{Solution to the optimal control problem}\label{sec:OT-solution}
The optimal control problem~\eqref{eq:opt-control} is precisely the Benamou-Brenier formulation of the optimal transport problem on a Riemannian manifold~\cite{benamou2000computational}. The optimal cost coincides with the optimal cost of the Kantorovich problem~\eqref{eq:Kantorovich}. 
By the Cauchy-Schwarz inequality we have
\begin{align*}
    \int_0^1 \|u(t)\|_g^2\,dt \geq (\int_0^1 \|u(t)\|_g\, dt)^2 = d(X(0),X(1))^2, 
\end{align*}
with equality when $\|u(t)\|_g$ is constant and the trajectory is a geodesic. Upon taking expectation of both sides and applying the initial and terminal constraints, we have that
\begin{align*}
    \mathbb E\left[\int_0^1 \|u(t)\|_g^2\,dt\right] \!\geq \!\mathbb E_{(X(0),X(1))\sim \pi}\left[d(X(0),X(1))^2\right]
\end{align*}
for any $\pi \in  \Pi(P,Q)$. Lastly, by infimizing both sides we obtain a lower-bound for the optimal control cost, 
with equality when $u(t)$ steers $X(t)$ along geodesics that connect $X(0)$ to $T(X(0))$ where $T(x)=\exp_{x}[-\nabla \phi(x)]$ is the optimal transport map from $P$ to $Q$. 
As a result, the solution to the optimal control problem is given by 
 \begin{align}\label{eq:optimal-trajectory}
     X(t) = \exp_{X(0)}[-t\nabla \phi(X(0))],\quad \text{for}\quad t\in[0,1]
 \end{align}
 where $\phi$ solves the dual Kantorovich problem~\eqref{eq:dual-Kantorovich}. The exponential exemplifies that the optimal trajectories are geodesics. 
 

In order to numerically solve the Kantorovich dual problem~\eqref{eq:dual-Kantorovich}, we use the result of Theorem~\ref{thm:McKann} to replace $\psi$ with $\phi^c$. 
Using the definition of $\phi^c$, 
\begin{align*}
 &\mathbb E_{Z'\sim Q}[\phi^c(Z')] = \mathbb E_{Z'\sim Q}[\min_z c(z,Z')-\phi(z)]\\
 & = \min_U \,\mathbb E_{Z'\sim Q}[ c(\exp_{Z'}[-U(Z')],Z')-\phi(\exp_{Z'}[-U(Z')])] 
\end{align*}
where  we assumed that  \[\exp_{Z'}[-U(Z')] = \arg \min_z c(z,Z') - \phi(z),\quad \text{a.e.}\]for some vector-field $U:\mathcal M \to T\mathcal M$. With this assumption, we can express the Kantorovich dual problem as
\begin{align}
\label{eq:max-min-OT}
&\max_\phi \,\min_U\, \mathbb E_{Z\sim P}[\phi(Z)] \\&+ \mathbb E_{Z'\sim Q}[ c(\exp_{Z'}(-U(Z')),Z')-\phi(\exp_{Z'}(-U(Z')))] \nonumber
\end{align}
This is a max-min optimization problem with the objective function that can be approximated using samples from the two distributions $P$ and $Q$. If $(\phi,U)$ is an optimal pair, then $T(z)=\exp_{z}[-\nabla \phi(z)]$ is the optimal transport map from $P$ to $Q$ and  $T^{-1}(z')=\exp_{z'}[-U(z')]$ is the optimal transport map from $Q$ to $P$.




 
\subsection{Solution to the optimal filtering problem}
The problem of computing conditional distributions can be formulated as an OT problem using the idea of block-triangular transport maps. In particular,  if a  map of the form $S:(x,y) \mapsto (T(x,y),y)$ transports the independent coupling $P_X \otimes P_Y$ to the joint distribution $P_{X Y}$, then the map $T(x,y)$ transports the prior $P_X$ to the conditional distribution $P_{X|Y}(\cdot|y)$ for any value $y$ of the observation; see Theorem 2.4 in~\cite{kovachki2020conditional} for a proof of this result. That is, 
\begin{align*}
    \text{if}&\quad (T,\text{Id})_{\#} (P_X \otimes P_Y) = P_{XY}\\
    \text{then}&\quad T(\cdot,y)_{\#}P_X = P_{X|Y}(\cdot|y),\quad \text{a.e.}~y.
\end{align*}
To seek the optimal transport map $T$ that characterizes posterior distributions, we express the Monge problem~\eqref{eq:Monge} with $P=P_X \otimes P_Y$, $Q=P_{XY}$, constrain the maps to have a block-triangular structure, form the Kantorovich dual, and apply the inf-$c$ convolution representation explained above in Section~\ref{sec:OT-solution}, which yields the max-min optimization problem
\begin{align}\label{eq:max-min-OT-FPF}
 \max_\phi \min_U\,&\mathbb E_{(X,Y)\sim P_{XY}}[\phi(X,Y)] \\&+ \mathbb E_{(X,Y) \sim P_X \otimes P_Y}[c(X,\exp_{X}(-T(X,Y)) \nonumber \\&- \phi(\exp_{X}(-T(X,Y))
,Y)]. \nonumber
\end{align}
Similarly to the optimization problem~\eqref{eq:max-min-OT}, the objective function in~\eqref{eq:max-min-OT-FPF} can be approximated using samples of the pair $(X,Y)\sim P_{XY}$. If $(\phi,U)$ is an optimal pair, then $T(x,y)=\exp_{z}[-U(x,y)]$ is the optimal transport map from $P_X$ to $P_{X|Y}(\cdot|y)$ for any value $y$ of the observation. 
\section{Numerical results}\label{sec:numerics}
In order to numerically solve the proposed optimal control~\eqref{eq:opt-control} and optimal filtering~\eqref{eq:Bayes} problems, we use their stochastic optimization formulations~\eqref{eq:max-min-OT} and~\eqref{eq:max-min-OT-FPF}, respectively. Both of these problems involve optimizing over a real-valued function $\phi:\mathcal M \times \mathbb R^m \to \mathbb R$, and a vector-field $U:\mathcal M\times \mathbb R^m \to  T\mathcal M$ (with the convention that $m=0$ for~\eqref{eq:max-min-OT}). To this end, we use neural networks to represent $\phi$ and $U$. The input layer to the neural network is based on a specific coordinate representation of the manifold that is described within each example. For all examples, the middle layers are fully connected residual blocks (for these examples, we used one or two blocks of size $32$). When $\mathcal M$ is a $n$-dimensional Lie-group, we model the output of the vector-field as an element of the Lie-algebra, which is isomorphic to $\mathbb R^n$.  We used the ADAM optimizer to solve the max-min problem with batch-size $64$, learning rate $10^{-3}$, $10$ inner-loop minimization iterations per one maximization iteration, and a total of $10^3$ to $10^4$ maximization iterations. The neural networks are trained on $10^3$ to $10^5$ samples from the distributions. The details of the numerical codes are available online\footnote{\url{https://github.com/Dan-Grange/OTManifold}}.

\subsection{Optimal transport mapping on $S^1$}
\begin{figure*}[t]
    \centering
    \subfigure[]{\includegraphics[width=0.22\textwidth,trim={100 100 100 100},clip]{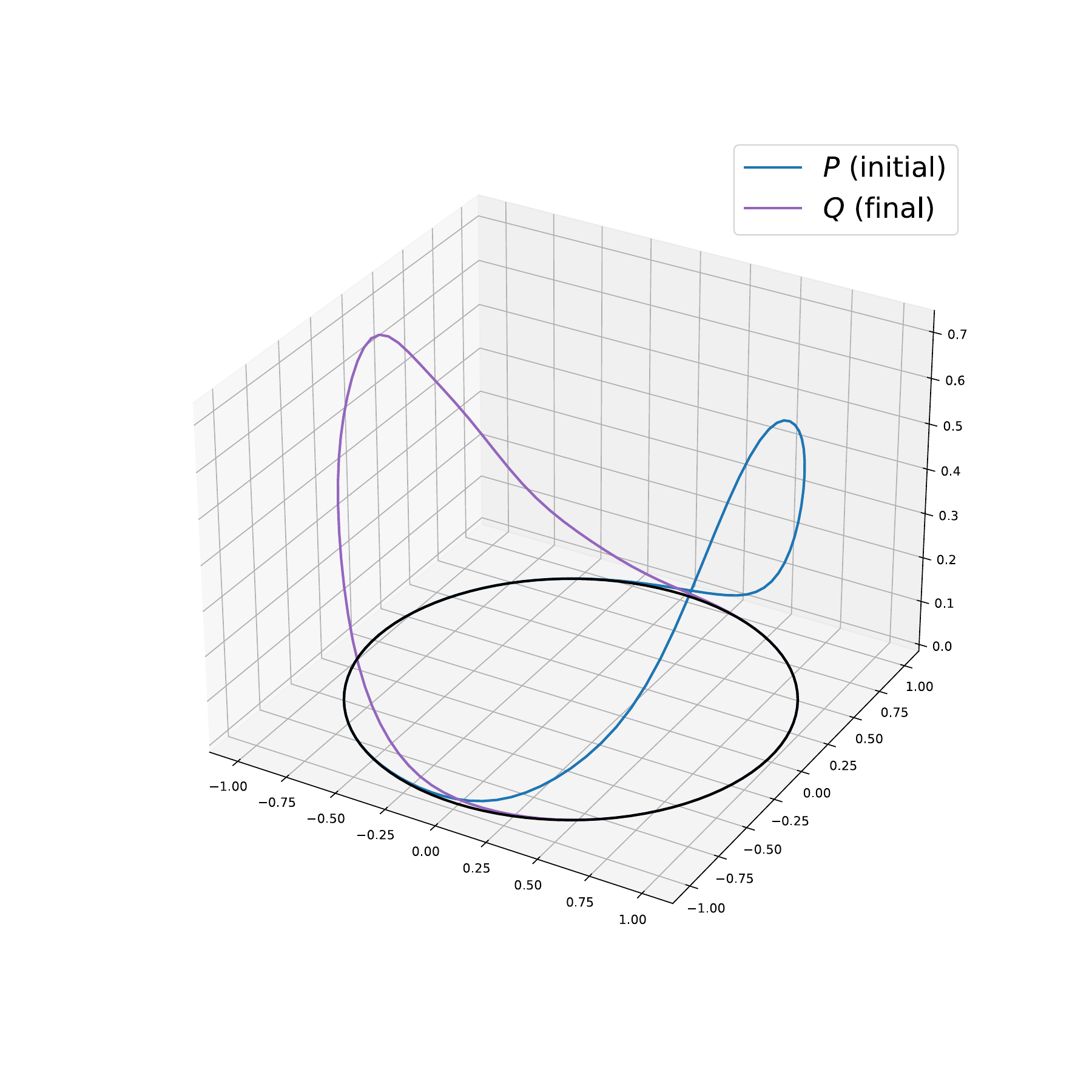}
     \includegraphics[width=0.22\textwidth,trim={100 100 100 100},clip]{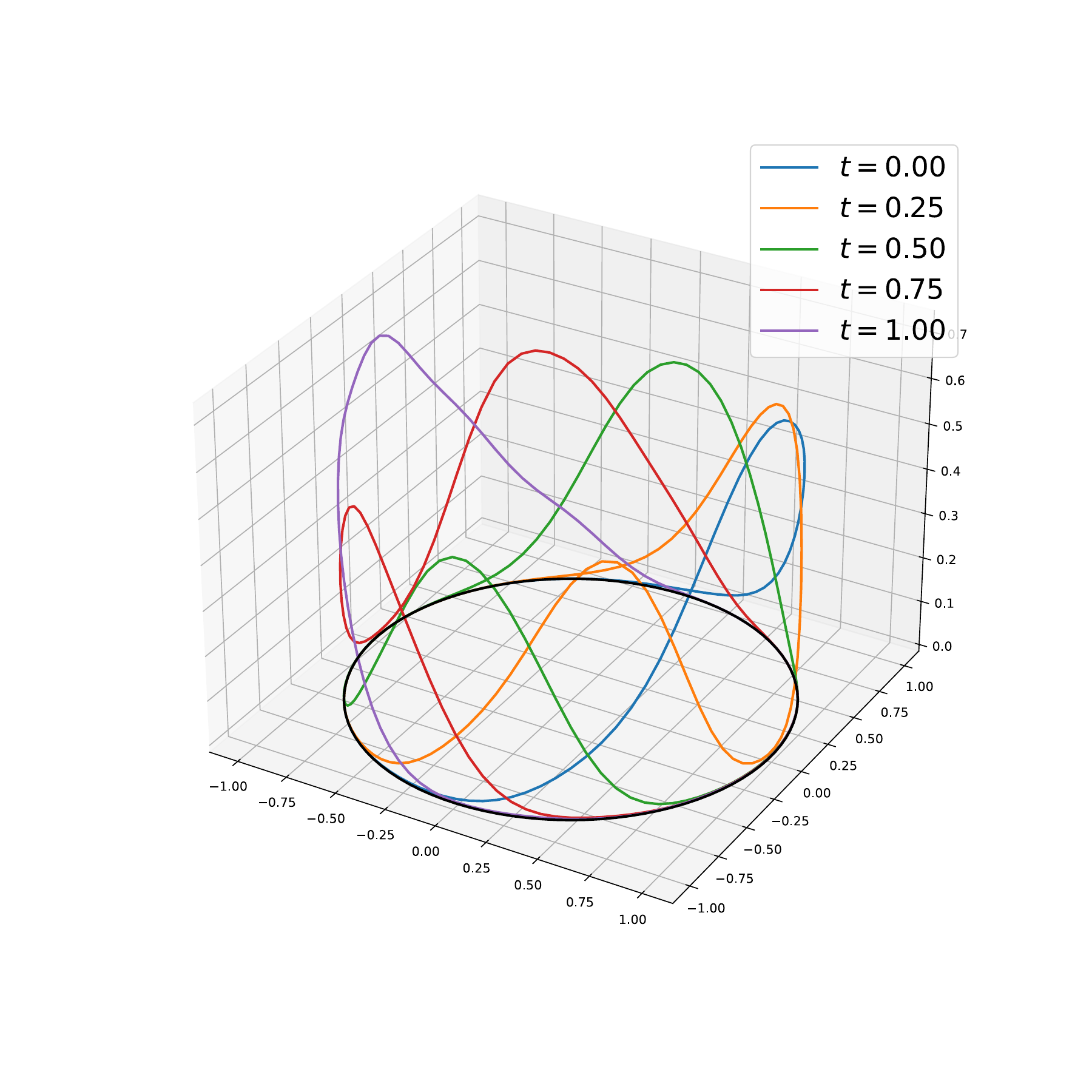}}\hfill 
     \subfigure[]{\includegraphics[width=0.22\textwidth,trim={100 100 100 100},clip]{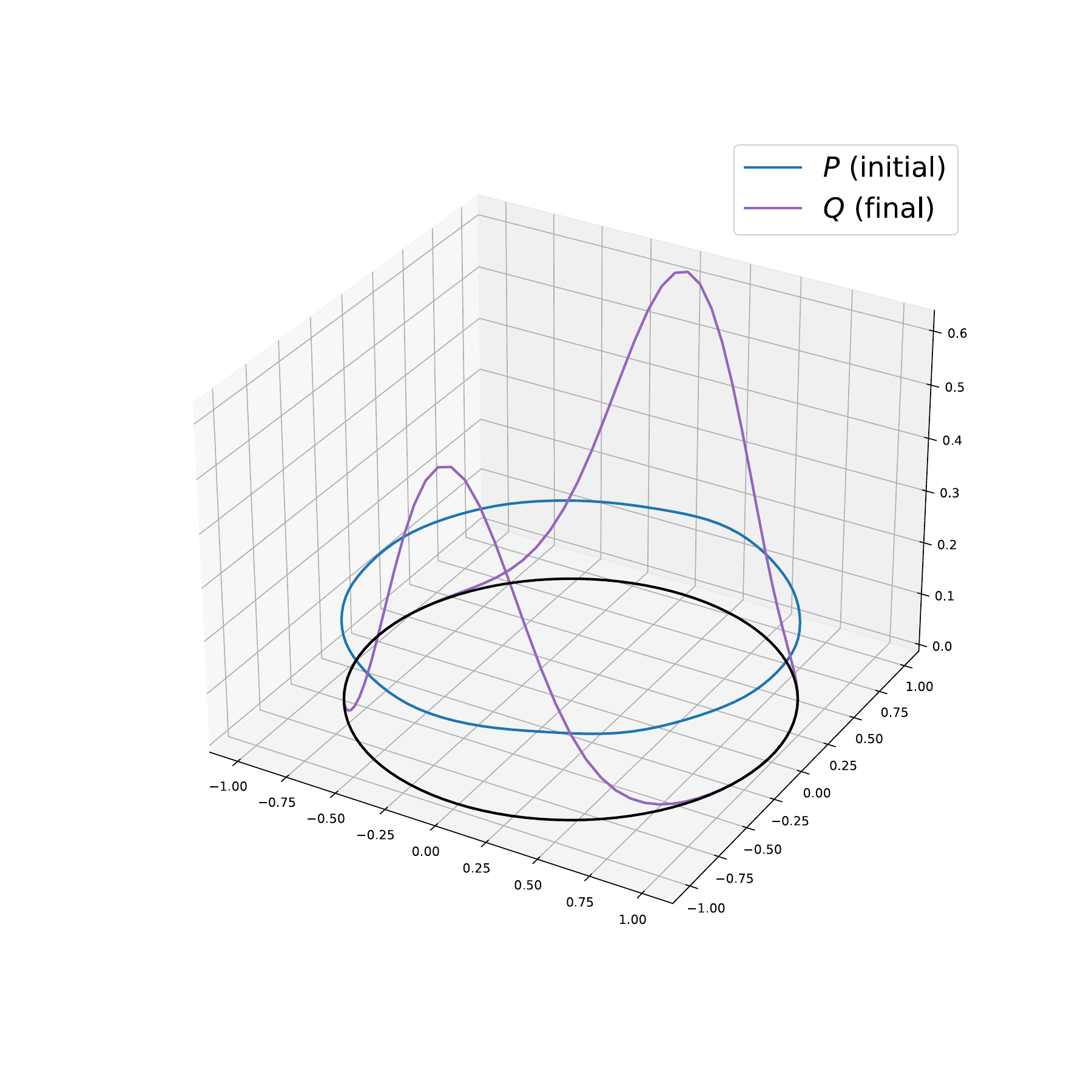}
    \includegraphics[width=0.22\textwidth,trim={100 100 100 100},clip]{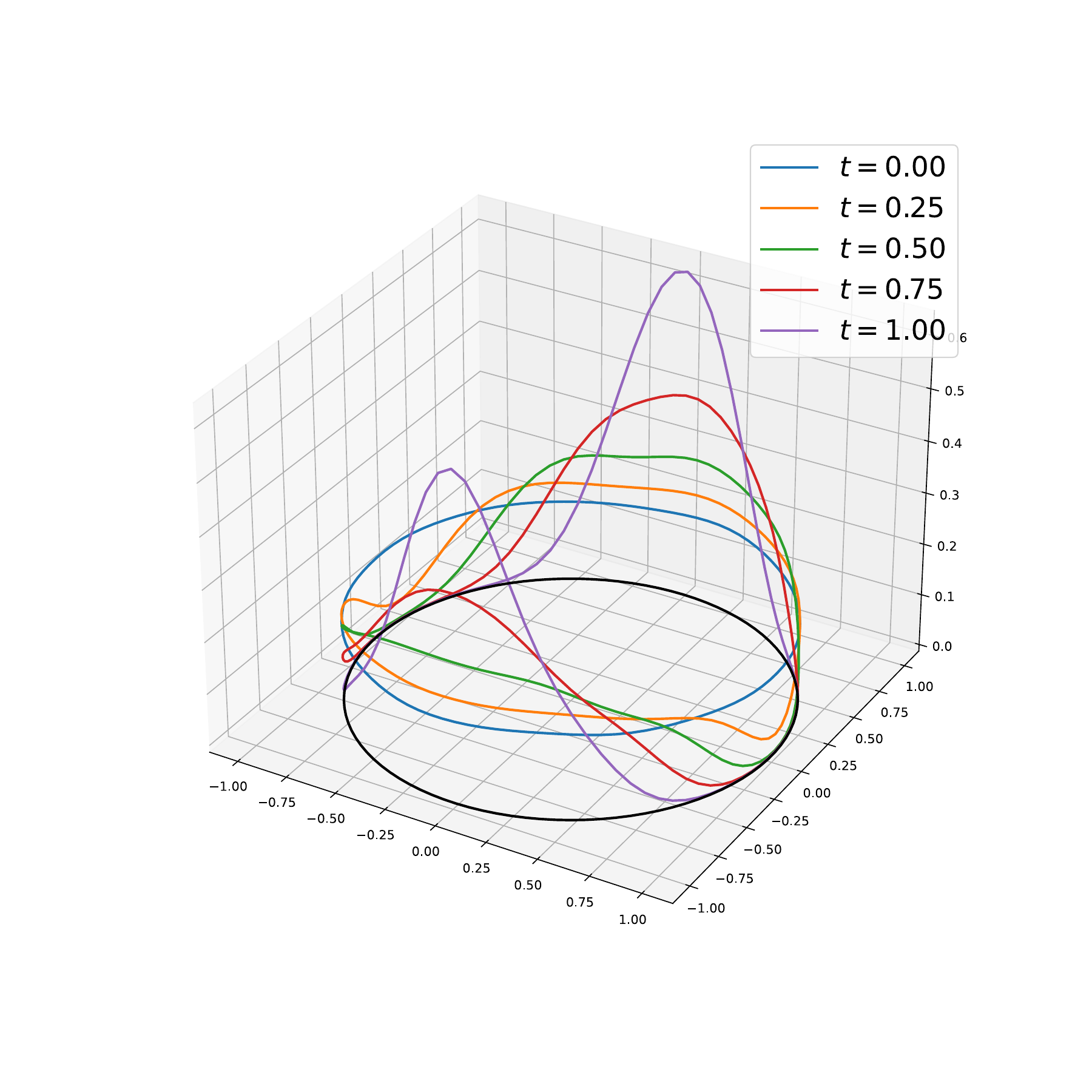}}
    \caption{OT mapping on the circle. The left panels in (a) and (b) show the given initial and final distributions. The right panels show the trajectory of optimal transportation evaluated by our numerical procedure. The curves are obtained with a kernel density estimation on samples.}
    \label{fig:OT-SO(2)}
\end{figure*}

The first example of OT on the circle $S^1$ uses coordinates $\theta \in [0,2\pi)$, with the tangent space $T_\theta S^1\cong \mathbb R^1$. For all $u \in T_\theta S^1$, we use the exponential map $\exp_x(u) = (x+u) \mod 2\pi$, and the geodesic distance 
\begin{align*}
d_{S^1}(\theta,\theta') =  \min(|\theta-\theta'|,2\pi-|\theta-\theta'|),\quad \forall \theta,\theta' \in [0,2\pi). 
\end{align*}
We implement our proposed numerical procedure for two choices of the marginal distribution: (1) $P$ and $Q$ are distributions centered at $\bar\theta=0$ and $\pi$, respectively, constructed from standard Gaussians modulo\footnote{We set $P,Q$ as having density $e^{-((\theta-\bar\theta) \,{\rm mod}2\pi)^2/2}/\sqrt{2\pi}$.} $2\pi$; (2) $P$ is the uniform distribution on $S^1$ and $Q$ is a mixture of two similarly projected Gaussians, centered at $0$ and $\pi$. We use the map $\theta \mapsto [\cos(\theta),\sin(\theta)]$ as the input layer of the networks. The given initial and final distribution, along with the optimal trajectory $\exp_x[-tU(x)]_{\#}P$ for $t\in[0,1]$ are depicted in Figure~\ref{fig:OT-SO(2)}. The first case highlights the periodic structure of the circle, where half of the initial distribution is transported in a clockwise direction, while the other half moves counter-clockwise, as opposed a constant shift by $\pi$ which is optimal in the Euclidean case. The second example highlights the ability of transporting to multi-modal distributions on the circle. 

\subsection{Optimal transport mapping on $SE(2)$}
\begin{figure}[t]
    \centering
    \includegraphics[width=0.24\textwidth,trim={50 50 50 50},clip]{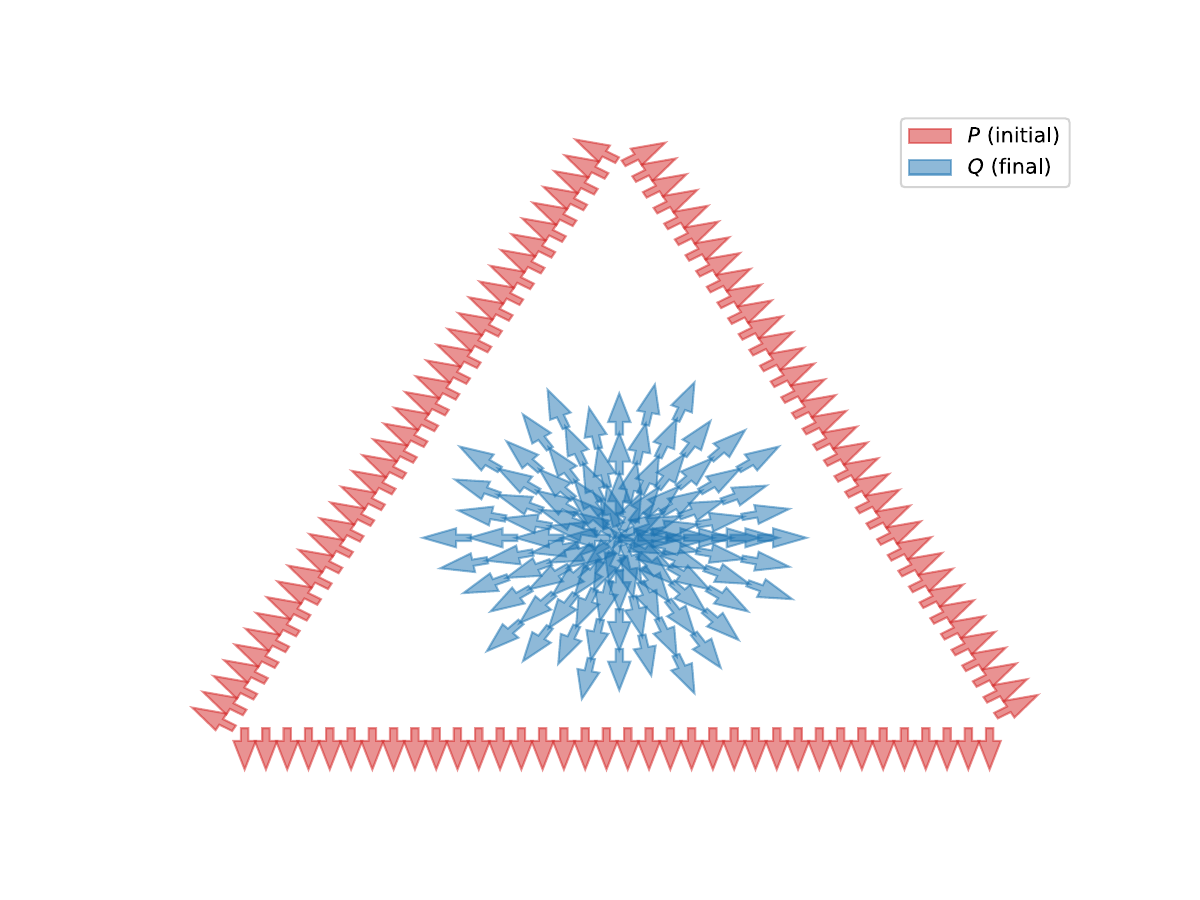}
     \includegraphics[width=0.24\textwidth,trim={50 50 50 50},clip]{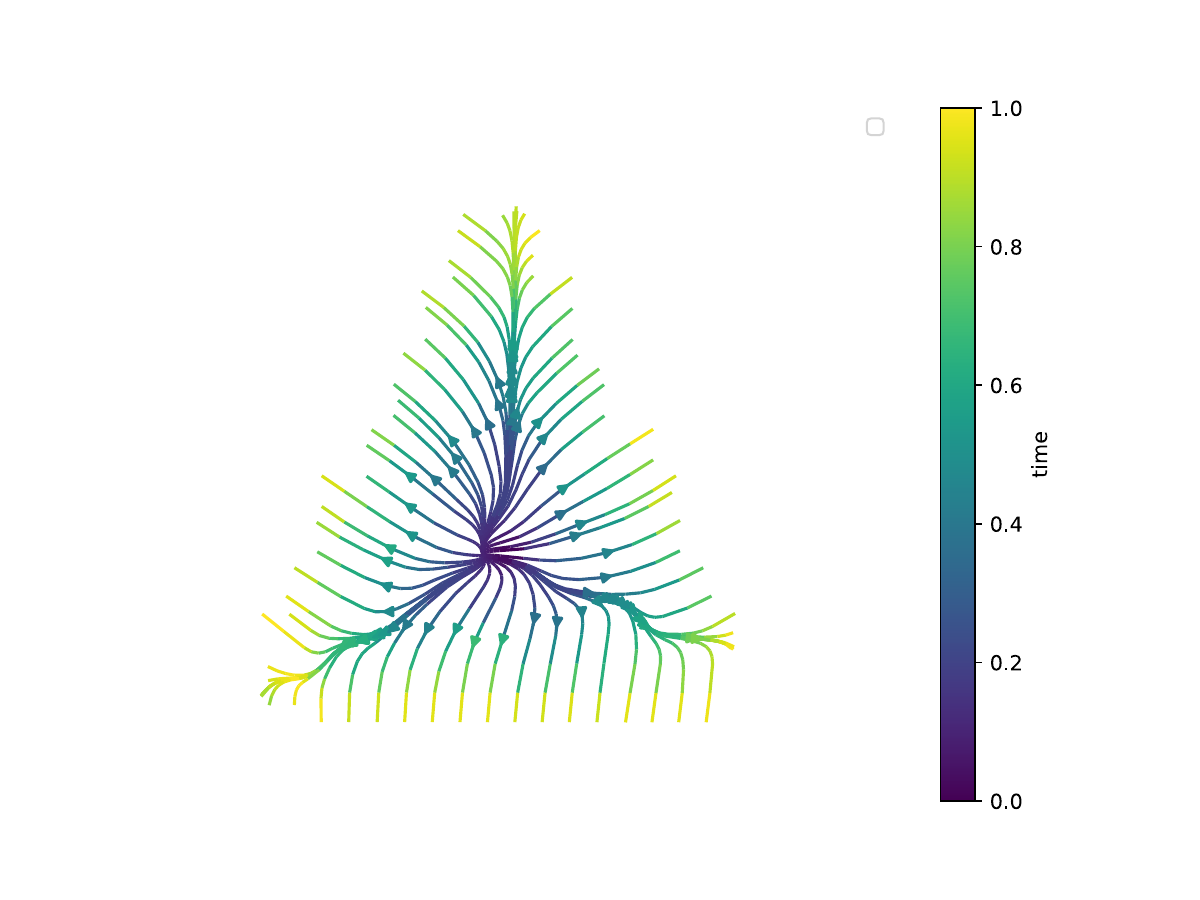}
    \caption{OT mapping on $SE(2)$. The elements of $SE(2)$ are depicted as arrows. Samples from the marginal distributions $P$ and $Q$ are shown in the left panel, and the optimal flow, obtained as a result of our numerical procedure, is depicted in the right panel. }
    \label{fig:OT-SE(2)}
\end{figure}

Next, we extend the previous example  to the special Euclidean group $SE(2)$. Elements of $SE(2)$ are represented by the coordinate $z=(x,y,\theta) \in \mathbb R^2 \times [0,2\pi)$. For all $u \in T_z SE(2)\cong \mathbb R^3$, we consider the exponential map
\[\exp_z(u) = (x+u_x,y+u_y,(\theta + u_\theta) \!\! \mod 2\pi),\]
and the geodesic distance 
\begin{align*}
d_{SE(2)}^2(z,z')^2 = d_{S^1}(\theta,\theta')^2 + |x - x'|^2 +  |y - y'|^2
\end{align*}
for $z=(x,y,\theta)$ and $z'=(x',y',\theta')$. 
We use $z \to [x,y,\cos(\theta),\sin(\theta)]$ as the input layer to the neural networks. Figure~\ref{fig:OT-SE(2)} qualitatively demonstrates the capability of our proposed approach to optimally transport distributions on $SE(2)$. 



\subsection{Optimal filtering on $SO(2)$} \label{sec:filtering_SO2}
\begin{figure}[h]
    \centering
    \includegraphics[width = .4\hsize,trim={190 60 210 90},clip]{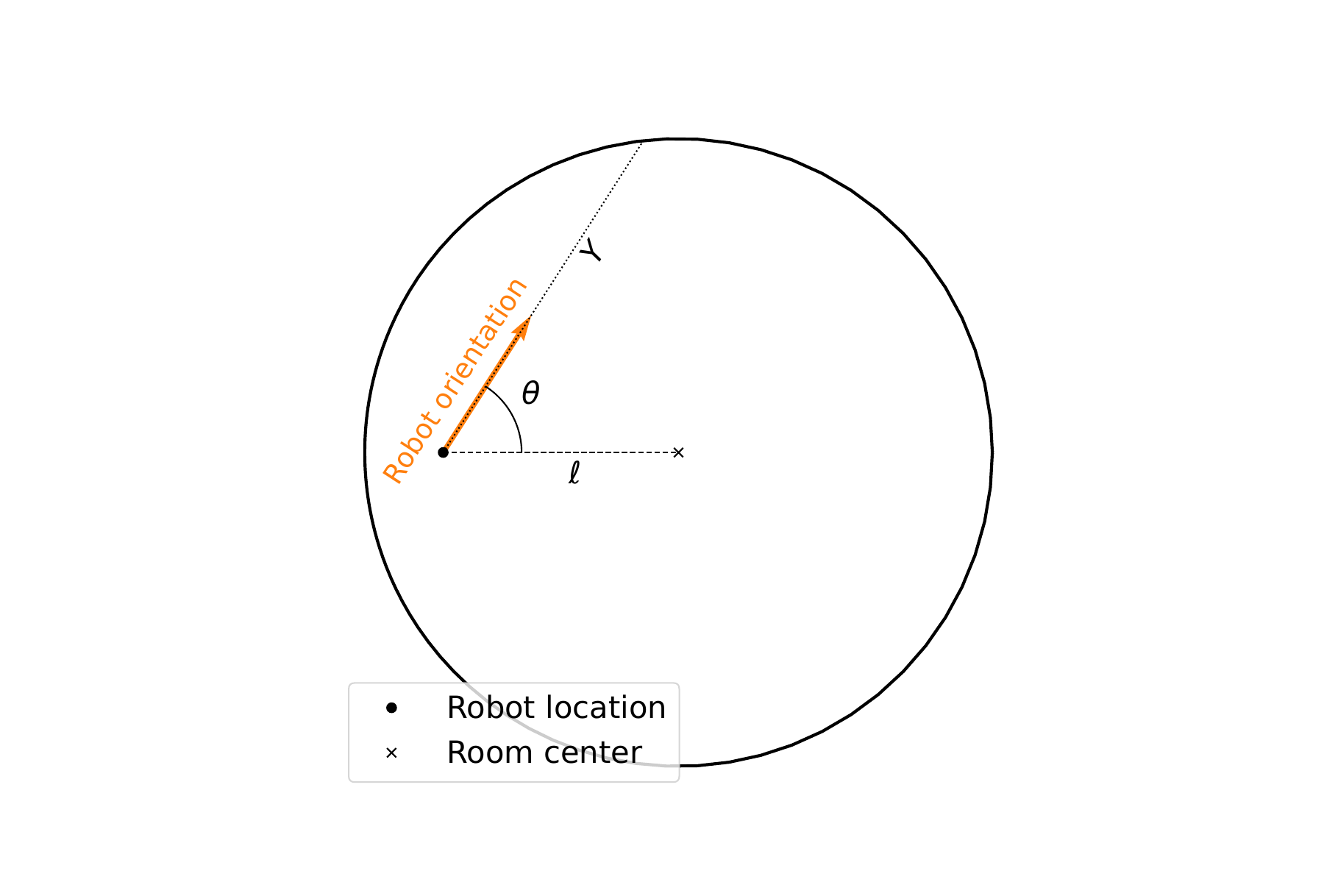}
    \includegraphics[width = .45\hsize,trim={220 60 180 80},clip]{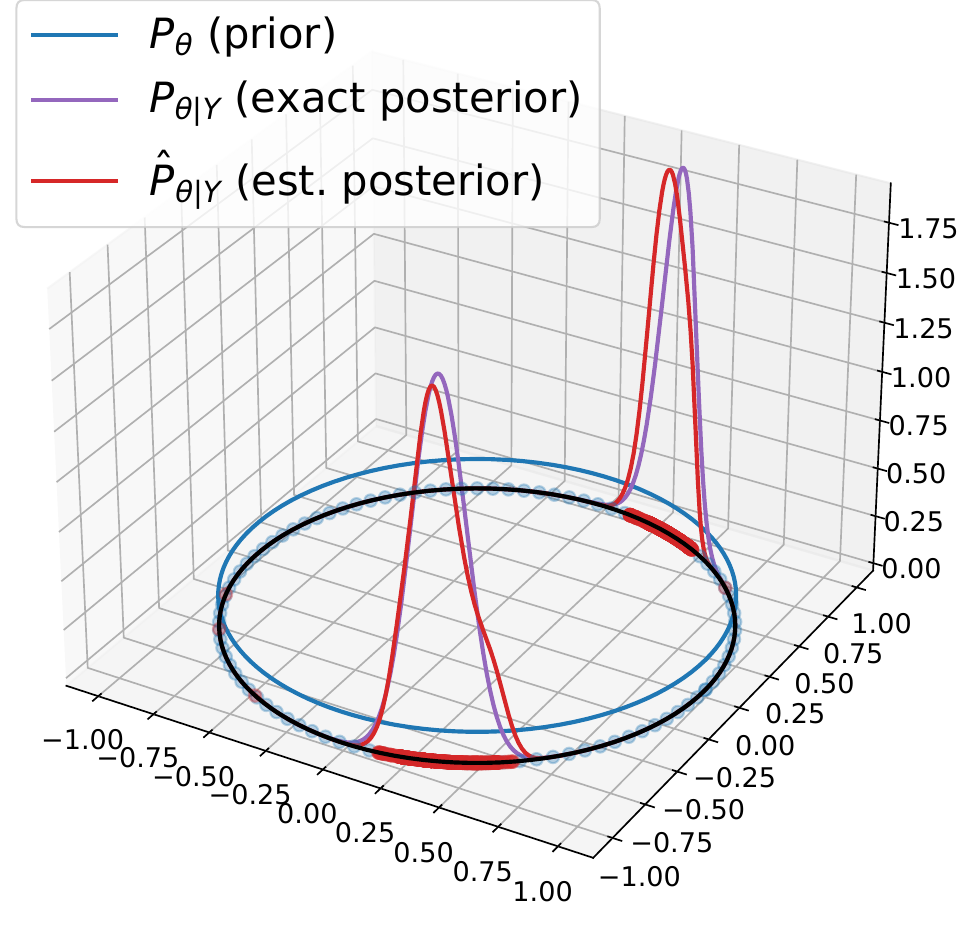}    
    \caption{Optimal filtering on a circle. The left panel shows the problem setup. The right panel shows the prior, exact and estimated posterior conditional distribution of the orientation given the observed signal. As expected, the bimodal posterior is captured by our algorithm. 
    }
    \label{fig:OT-FPF-SO2-static}
\end{figure}

We now consider the problem of estimating the orientation $\theta\in [0,2\pi)$ of a ground robot in a circular room as depicted in Figure~\ref{fig:OT-FPF-SO2-static}. The robot is equipped with a sensor that measures its distance to the wall of the room along the direction that the robot is facing, i.e. $Y=h(\theta) + W$, where 
\begin{align*}
    h(\theta) = \ell \cos(\theta) + \sqrt{1-\ell^2\sin(\theta)^2}, 
\end{align*}
$\ell$ is the distance of the robot from the center, and $W\sim N(0,10^{-2})$ is additive Gaussian noise. We assume a uniform prior distribution for the orientation $\theta$ and use our proposed numerical procedure to characterize the conditional distribution of the orientation $\theta$ given the observation $Y$.   
We use the input layer $(\theta,y) \mapsto (\cos(\theta),\sin(\theta),y)$ in our networks for this example. Figure~\ref{fig:OT-FPF-SO2-static} depicts the resulting posterior distribution. Due to the symmetric geometry of the circle, there are two angles that are consistent with an observation. As a result, the conditional distribution is  bimodal, which is approximated well by our numerical procedure.   

\begin{figure}[t]
\centering
\includegraphics[width = .5\hsize,trim={0 0 30 15},clip]{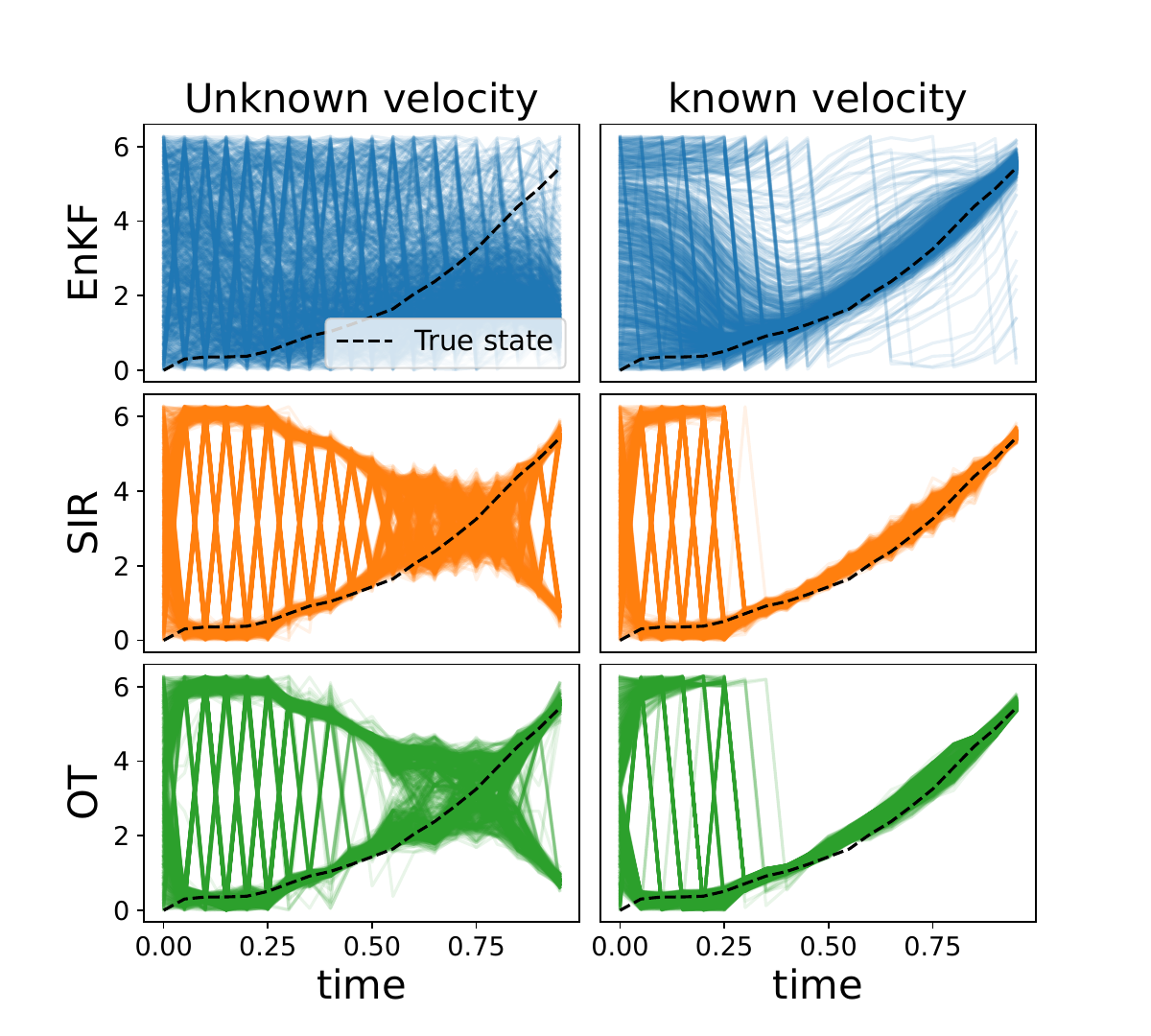}
\includegraphics[width = .45\hsize,trim={0 0 0 0},clip]{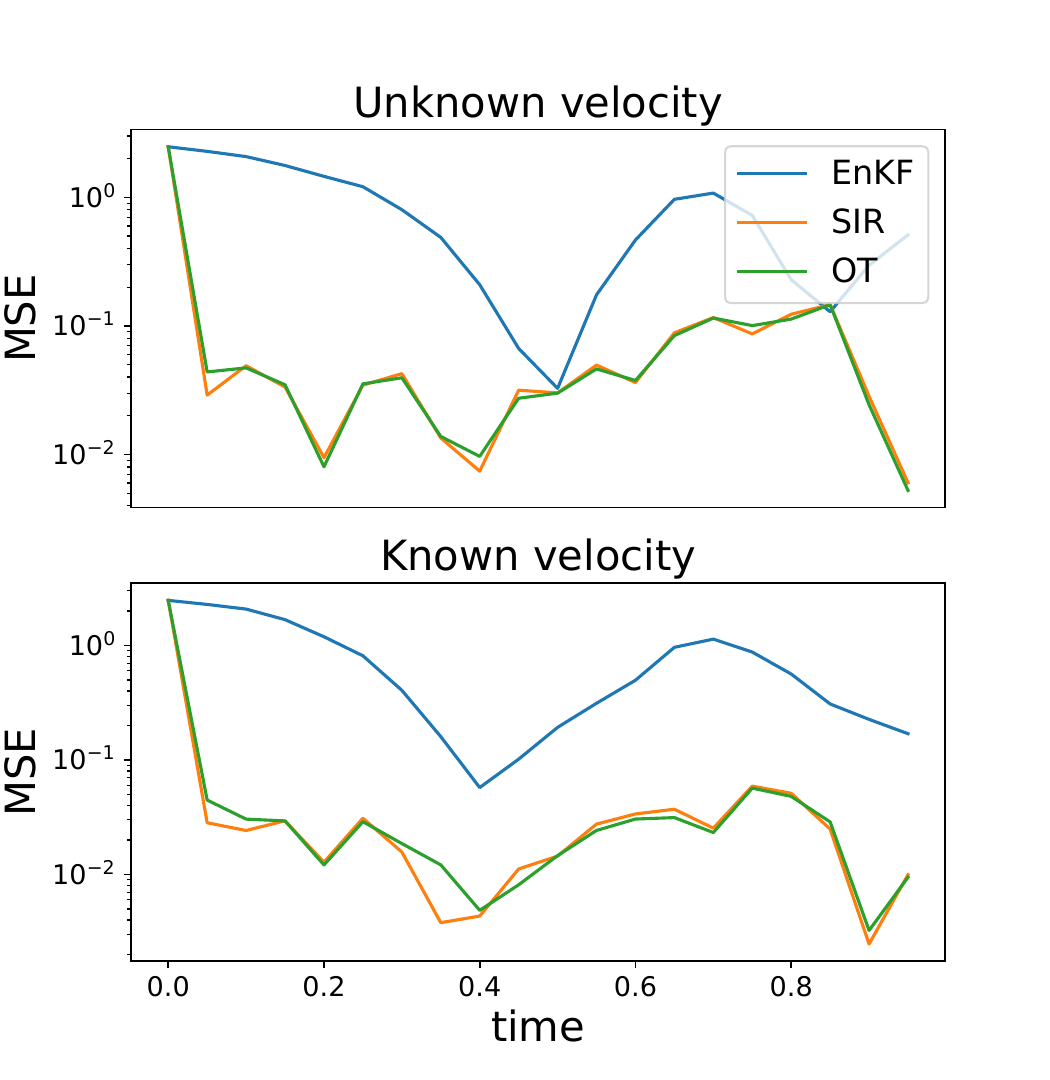}
\caption{Simulation results of EnKF, SIR, and OT for the (dynamic) optimal filtering example on $SO(2)$. Left panel: the trajectory of the samples/particles along with the true state for the cases of known and unknown velocity. Right panel: The mean-squared-error (MSE) for estimating the true state averaged over 10 independent simulations.}
\label{fig:OT-FPF-SO2-dynamic}
\end{figure}

In order to test the algorithm in a dynamic setting, we let the orientation change with a constant velocity $u$ using the model: $\theta_{t+1}=\theta_t + u + \xi_t$ where $\xi_t \sim N(0,10^{-2})$ is the process noise. Our goal is now to compute the conditional (filtering) distribution of the orientation $\theta_t$ given the history of sensory observation $Y_1,Y_2,\ldots,Y_t$, generated by $Y_t=h(\theta_t)+W_t$. To do so, we use the numerical framework in Section~\ref{sec:methodology} to update the conditional distributions for $\theta_t$ sequentially as new observations arrive, as in~\cite{al2023optimal}.  
We consider two setups for the orientation's  velocity $u$ in this example. In the first setup, we assume the velocity is known and non-zero. In this case, due to the known direction of the motion, the problem becomes observable, leading to a unimodal filtering distribution. In the second setup, we assume zero-velocity in our algorithm, while the actual velocity is non-zero. Since the direction of motion is unknown, the conditional distribution remains bimodal. The results are depicted in Figure~\ref{fig:OT-FPF-SO2-dynamic}. For comparison, we have shown results from the Ensemble Kalman filter (EnKF)~\cite{evensen2006} and the sequential importance sampling and resampling (SIR) particle filter~\cite{doucet09}. As a  quantitative comparison, we evaluated the mean-squared-error (MSE) in estimating  
$ f(\theta)=d_{S^1}(\theta,0).$
The result shows consistency of our algorithm with SIR. We expect the OT approach to outperform SIR in high-dimensional problems, as shown in the Euclidean case in~\cite{al2023optimal}.

\subsection{Optimal filtering on $SE(2)$}
\begin{figure}[h]
    \centering
    \includegraphics[width = .4\hsize,trim={190 60 180 80},clip]{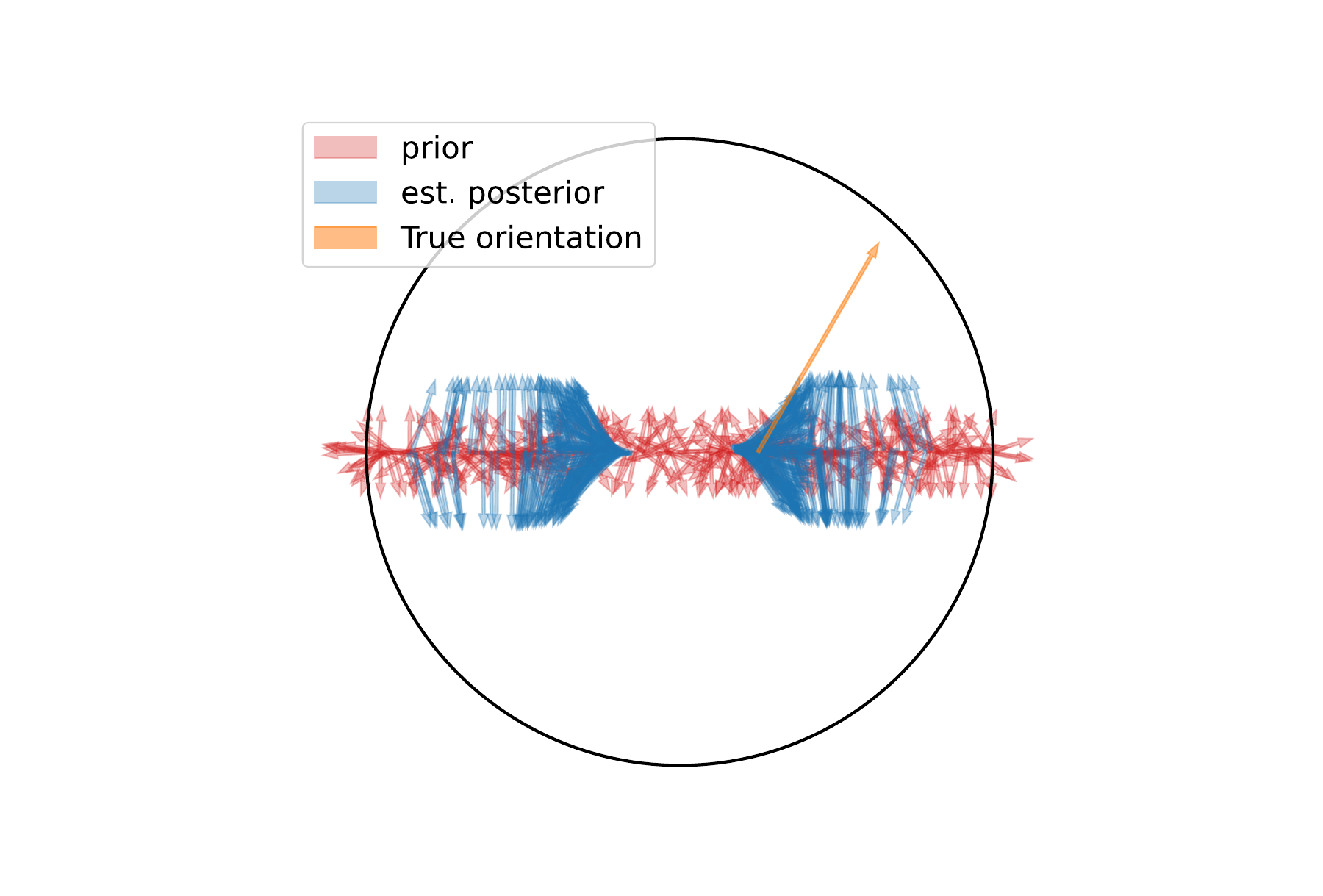}
    \includegraphics[width = .5\hsize,trim={50 20 80 20},clip]{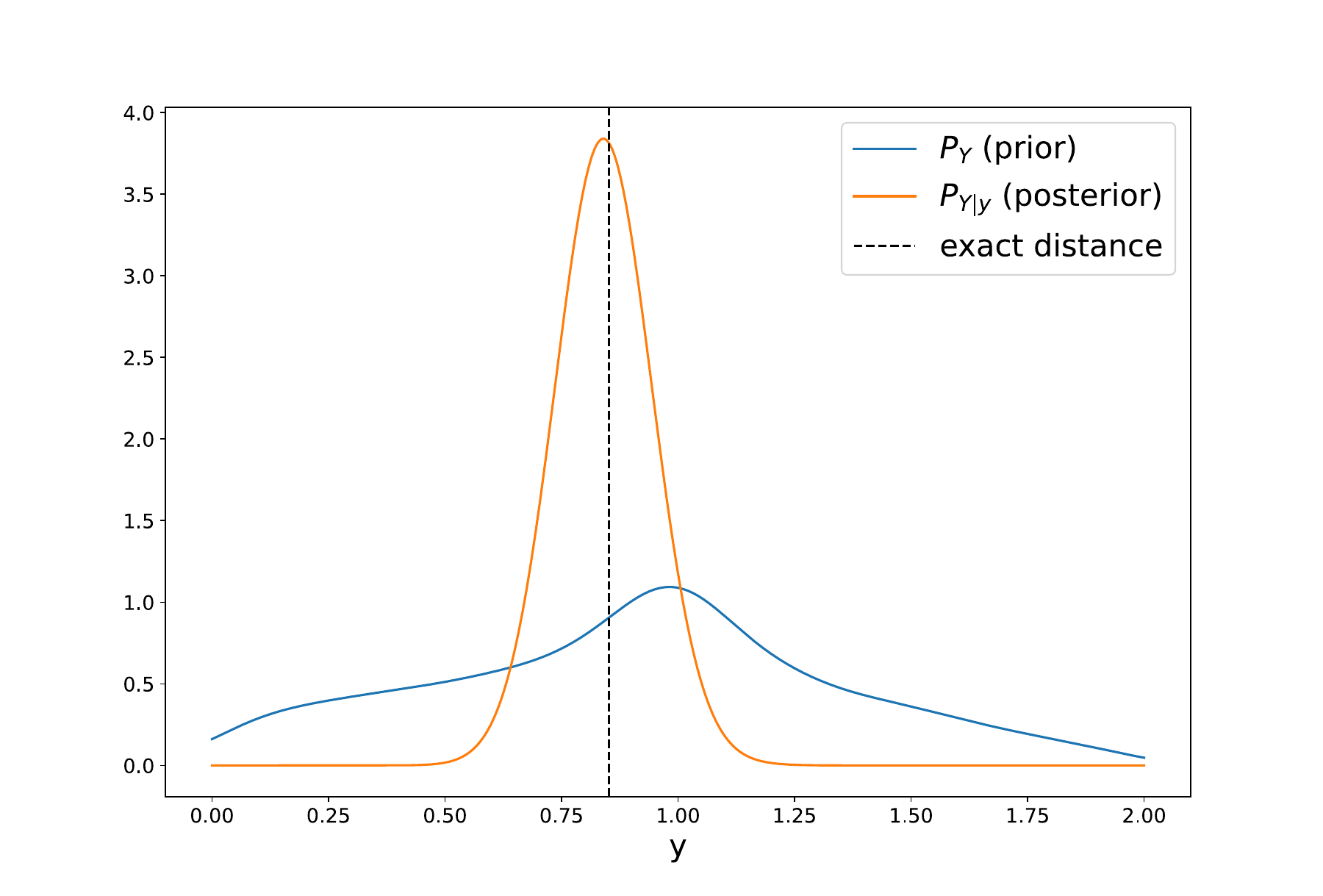}    
    \caption{Optimal filtering on SE(2). The left panel shows samples from the prior and posterior distribution of the location and orientation of the robot. The right panel shows the prior and posterior prediction for the observation signal, in comparison to the value of the exact distance. 
    }
    \label{fig:OT-FPF-SE2-static}
\end{figure}

Building on the previous example, we assume the position of the robot is also unknown, and consider the goal of estimating the position and orientation of the robot simultaneously. 
The state of the robot is now represented by $z=(x,y,\theta)$, an element of $SE(2)$. We assume $y=0$, and consider independent uniform distributions on $[-1,1]$ for $x$ and $[0,2\pi)$ for $\theta$. We consider the same distance measuring sensor as in Section~\ref{sec:filtering_SO2}. Figure \ref{fig:OT-FPF-SE2-static} depicts the simulation results. We observe that the algorithm is able to capture the complicated posterior distribution. For verification, we show the predictive distribution for the observation $Y$ prior to receiving the realized observation, and after receiving the observation. We note that the posterior predictive distribution is concentrated around the exact distance.

\subsection{Optimal transport mapping on $SO(3)$}
\begin{figure}[t]
    \centering
    \includegraphics[width=0.45\textwidth,trim={50 20 50 20},clip]{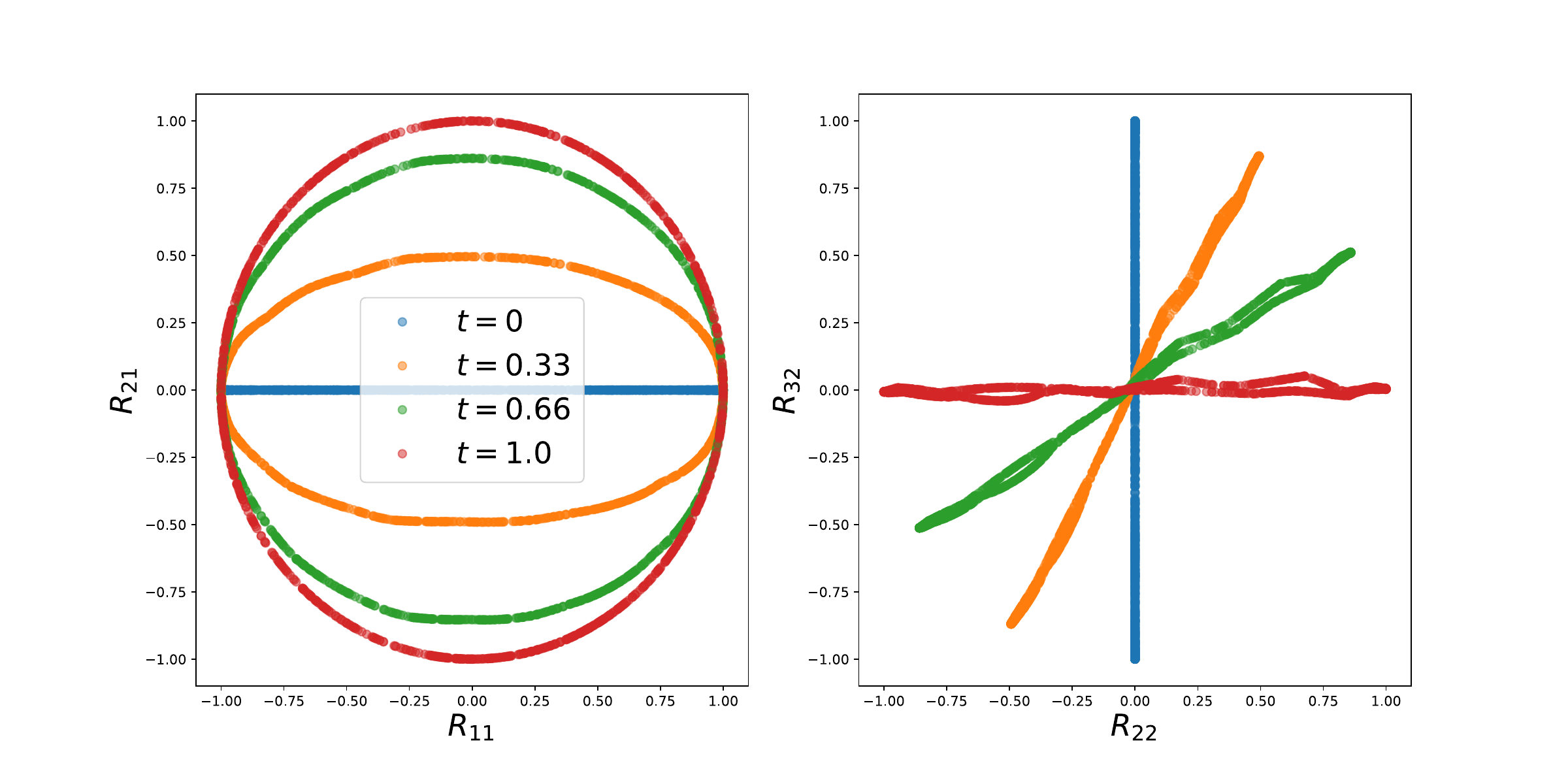}
    \caption{OT mapping on $SO(3)$.  Samples from the optimal transport trajectory from $P$ (blue) to $Q$ (red) are shown in the figure. The two figures show four components of the rotation matrix.}
    \label{fig:OT-SO(3)}
\end{figure}

Next, we consider the OT problem on the special orthogonal group $SO(3)$. An element of $SO(3)$ is represented by a rotation matrix $R$. The elements of the tangent space $T_RSO(3)$ are represented by $R\Omega$ where $\Omega$ is a skew-symmetric matrix. 
A three dimensional vector $\omega \in \mathbb R^3$ is uniquely mapped to a skew-symmetric matrix according to
\[
[\omega]_\times := \begin{bmatrix}
     \phantom{-}0 &-\omega_3 & \phantom{-}\omega_2\\
     \phantom{-}\omega_3 & \phantom{-}0 & -\omega_1 \\
     -\omega_2 & \phantom{-}\omega_1 & \phantom{-}0
\end{bmatrix}.
\]
We have the exponential map $\exp_R(R[\omega]_\times) = Re^{[\omega]_\times}$, where \[e^{[\omega]_\times} = I + \frac{\sin(\|\omega\|)}{\|\omega\|}[\omega]_\times +   \frac{1-\cos(\|\omega\|)}{\|\omega\|^2}[\omega]_\times^2.\]
Then, the geodesic distance is given by
\begin{align*}
    d_{SO(3)}(R,\exp_R(R[\omega]_\times))= \sqrt{2} \|\omega\|.
\end{align*}
The input to the neural net is the matrix representation of the rotation matrix. The output of the vector-field is designed to be a three dimensional vector $\omega$, which is then mapped to an element of the tangent space by $R[\omega]_\times$. Figure~\ref{fig:OT-SO(3)} depicts the result for optimally transporting an initial distribution of rotation matrices $P$ to a final distribution $Q$. The distributions $P$ and $Q$ are chosen to be
\begin{align*}
    \begin{bmatrix}
        &\cos U &\sin U &0\\
        &0           &0            &1\\
        &\sin U &-\cos U  &0\\
    \end{bmatrix},
    \begin{bmatrix}
        &-\cos U &-\sin U &0\\
        &-\sin U &\cos U  &0\\
        &0           &0            &-1\\
    \end{bmatrix}, \quad 
\end{align*}
respectively, with $U$ being uniformly distributed on $[0, 2\pi]$.
\begin{figure}[t]
\centering
\subfigure{\includegraphics[width = .32\hsize, trim={20 20 20 20},clip]{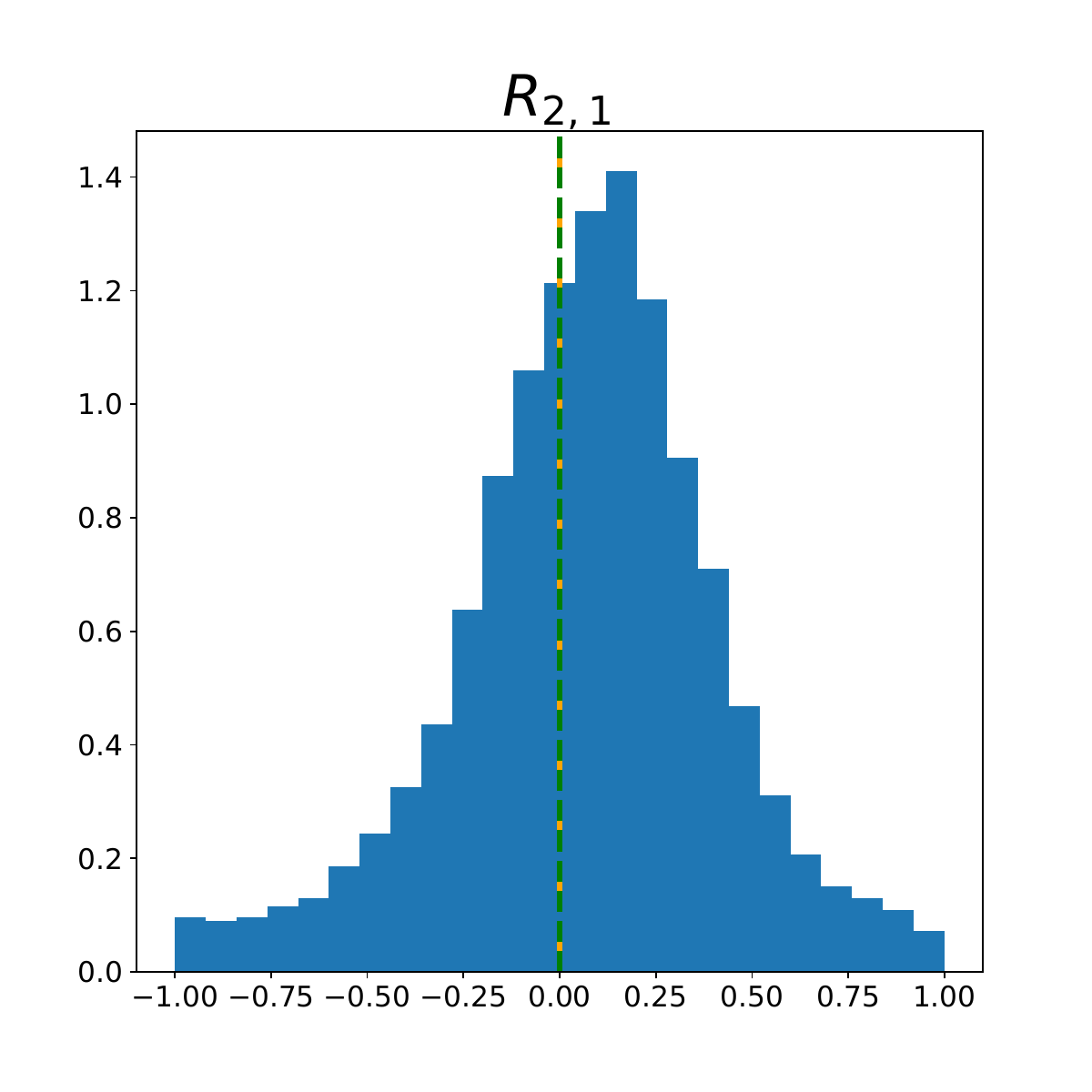}}\hfill 
\subfigure{\includegraphics[width = .32\hsize, trim={20 20 20 20},clip]{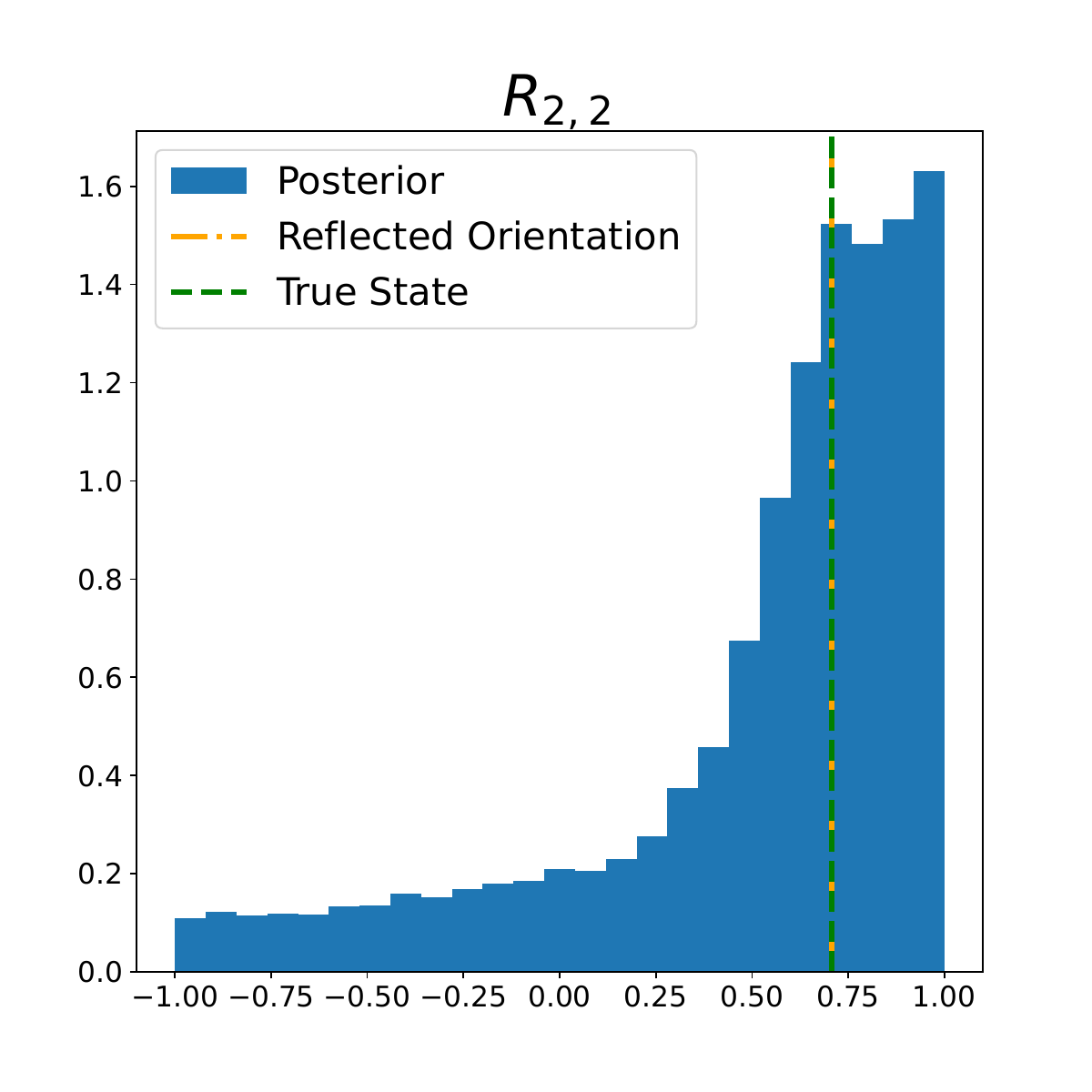}}\hfill 
\subfigure{\includegraphics[width = .32\hsize, trim={20 20 20 20},clip]{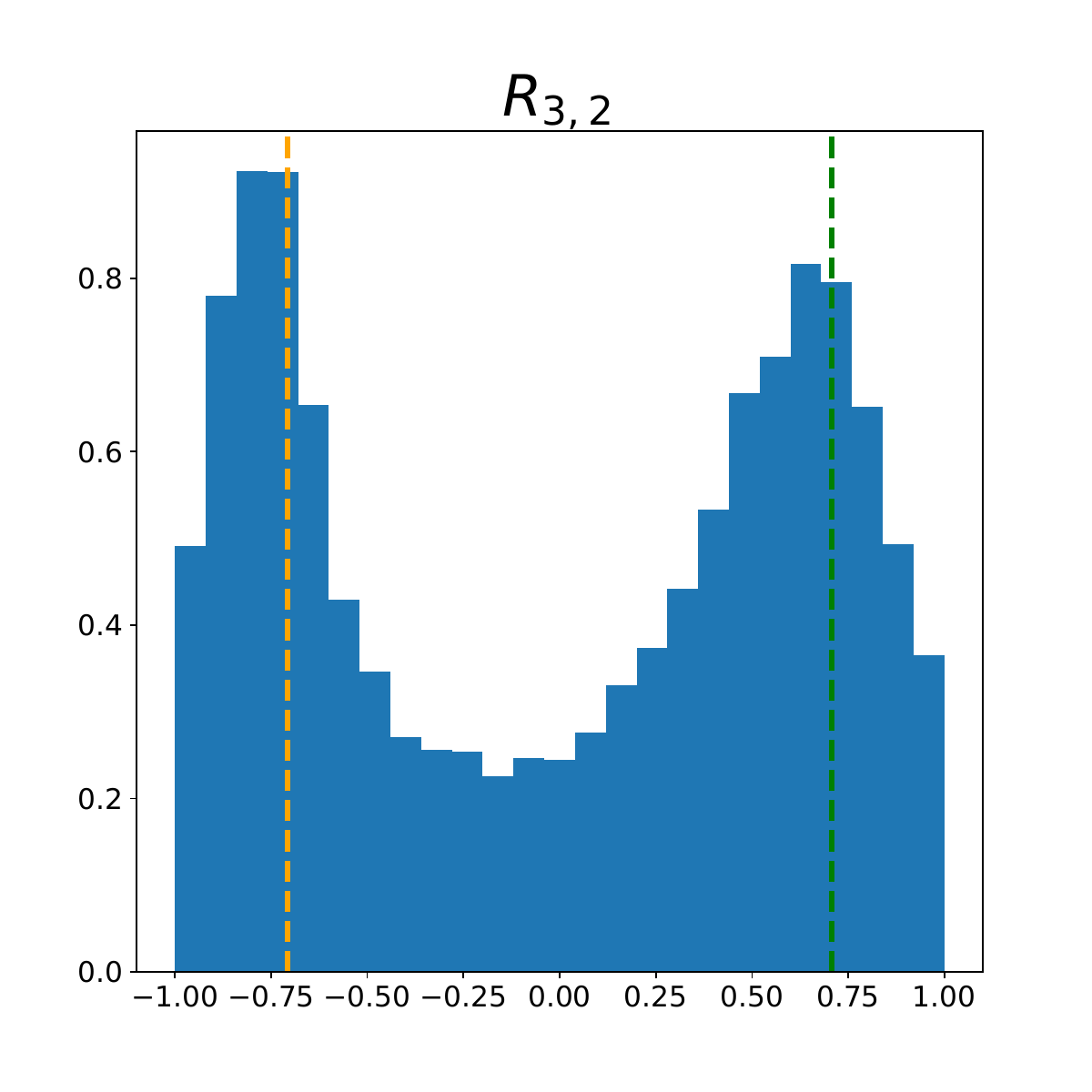}}
 \caption{Optimal filtering on $SO(3)$. The histograms show the approximate posterior distributions for three components of $R$ using samples computed by the estimated transport map. The dashed lines correspond to two rotation matrices that are consistent with the measurement.}
\label{fig:OT-FPF-SO(3)}
\end{figure}
\subsection{Optimal filtering on $SO(3)$}
Lastly, we consider the problem of characterizing the conditional distribution on $SO(3)$. 
We consider noisy measurement $Y=h(R) + W$ where 
$
    h(R) = (R_{1,1}, R_{2,1}, R_{3,1}, R_{2,2}).
$
and $W \sim N(0,10^{-2})$. The selected observation function is degenerate because there are two rotation matrices that correspond to any exact measurement. For example,  
\begin{align*}
    R=\begin{bmatrix}
        1 & 0 & 0\\
        0 & \frac{1}{\sqrt{2}} & \mp\frac{1}{\sqrt{2}} \\
        0 & \pm\frac{1}{\sqrt{2}} & \frac{1}{\sqrt{2}}
    \end{bmatrix}
\end{align*}
 are both consistent with $h(R)=(1,0,0,\frac{1}{\sqrt{2}})$. In our experiment, we consider a uniform prior probability distribution on  $SO(3)$ and use our proposed numerical procedure to sample from the conditional distribution. The numerical result is depicted in Figure \ref{fig:OT-FPF-SO(3)}. It is observed that the algorithm is able to capture the bimodal distribution.  

\section{Concluding remarks}
This work proposes a computational methodology for finding optimal transport maps and sampling conditional probability distributions on several instances of Riemannian manifolds. Our planned future research will focus on: (1) theoretical stability and sample-complexity analysis of the proposed approach, and (2) application of the approach to problems of attitude estimation using the publicly available Multi-Vehicle Stereo Event Camera (MVSEC) dataset~\cite{Zhu2018}. 

\bibliographystyle{IEEEtran}
\bibliography{refs}

\begin{thebibliography}{10}
\providecommand{\url}[1]{#1}
\csname url@samestyle\endcsname
\providecommand{\newblock}{\relax}
\providecommand{\bibinfo}[2]{#2}
\providecommand{\BIBentrySTDinterwordspacing}{\spaceskip=0pt\relax}
\providecommand{\BIBentryALTinterwordstretchfactor}{4}
\providecommand{\BIBentryALTinterwordspacing}{\spaceskip=\fontdimen2\font plus
\BIBentryALTinterwordstretchfactor\fontdimen3\font minus \fontdimen4\font\relax}
\providecommand{\BIBforeignlanguage}[2]{{%
\expandafter\ifx\csname l@#1\endcsname\relax
\typeout{** WARNING: IEEEtran.bst: No hyphenation pattern has been}%
\typeout{** loaded for the language `#1'. Using the pattern for}%
\typeout{** the default language instead.}%
\else
\language=\csname l@#1\endcsname
\fi
#2}}
\providecommand{\BIBdecl}{\relax}
\BIBdecl

\bibitem{al2023optimal}
M.~Al-Jarrah, B.~Hosseini, and A.~Taghvaei, ``Optimal transport particle filters,'' \emph{arXiv:2304.00392}, 2023.

\bibitem{sepulchre2021optimal}
R.~Sepulchre, ``Optimal transport in estimation and control [about this issue],'' \emph{IEEE Control Systems Magazine}, vol.~41, no.~4, pp. 6--8, 2021.

\bibitem{Vil08}
C.~Villani, \emph{Optimal {T}ransport: {O}ld and {N}ew}.\hskip 1em plus 0.5em minus 0.4em\relax Springer, 2008, vol. 338.

\bibitem{peyre2019computational}
G.~Peyr{\'e}, M.~Cuturi \emph{et~al.}, ``Computational optimal transport: With applications to data science,'' \emph{Foundations and Trends{\textregistered} in Machine Learning}, vol.~11, no. 5-6, pp. 355--607, 2019.

\bibitem{arjovsky2017wasserstein}
M.~Arjovsky, S.~Chintala, and L.~Bottou, ``{W}asserstein generative adversarial networks,'' in \emph{International conference on machine learning}.\hskip 1em plus 0.5em minus 0.4em\relax PMLR, 2017, pp. 214--223.

\bibitem{courty2015optimal}
N.~Courty, R.~Flamary, D.~Tuia, and A.~Rakotomamonjy, ``Optimal transport for domain adaptation,'' \emph{arXiv:1507.00504}, 2015.

\bibitem{kolouri2017optimal}
S.~Kolouri, S.~R. Park, M.~Thorpe, D.~Slepcev, and G.~K. Rohde, ``Optimal mass transport: Signal processing and machine-learning applications,'' \emph{IEEE {S}ignal {P}roc.\ {M}agazine}, vol.~34, no.~4, pp. 43--59, 2017.

\bibitem{dominitz2010texture}
A.~Dominitz and A.~Tannenbaum, ``Texture mapping via optimal mass transport,'' \emph{IEEE transactions on visualization and computer graphics}, vol.~16, no.~3, pp. 419--433, 2010.

\bibitem{rabin2011wasserstein}
J.~Rabin, G.~Peyr{\'e}, J.~Delon, and M.~Bernot, ``{W}asserstein barycenter and its application to texture mixing,'' in \emph{International Conference on Scale Space and Variational Methods in Computer Vision}.\hskip 1em plus 0.5em minus 0.4em\relax Springer, 2011, pp. 435--446.

\bibitem{su2015optimal}
Z.~Su, Y.~Wang, R.~Shi, W.~Zeng, J.~Sun, F.~Luo, and X.~Gu, ``Optimal mass transport for shape matching and comparison,'' \emph{IEEE transactions on pattern analysis and machine intelligence}, vol.~37, no.~11, pp. 2246--2259, 2015.

\bibitem{chen2021optimal}
Y.~Chen, J.~Karlsson, and A.~Ringh, ``Optimal transport for applications in control and estimation,'' \emph{IEEE Control Systems Magazine}, vol.~41, no.~4, pp. 28--33, 2021.

\bibitem{chen2021optimal2}
Y.~Chen, T.~T. Georgiou, and M.~Pavon, ``Optimal transport in systems and control,'' \emph{Annual Review of Control, Robotics, and Autonomous Systems}, vol.~4, pp. 89--113, 2021.

\bibitem{taghvaei2021OTFPF}
A.~Taghvaei and P.~G. Mehta, ``Optimal transportation methods in nonlinear filtering,'' \emph{IEEE Control Systems Magazine}, vol.~41, no.~4, pp. 34--49, 2021.

\bibitem{CheGeoTan18}
Y.~Chen, T.~T. Georgiou, and A.~Tannenbaum, ``Optimal transport for gaussian mixture models,'' \emph{IEEE Access}, vol.~7, pp. 6269--6278, 2018.

\bibitem{leygonie2019adversarial}
J.~Leygonie, J.~She, A.~Almahairi, S.~Rajeswar, and A.~Courville, ``Adversarial computation of optimal transport maps,'' \emph{arXiv:1906.09691}, 2019.

\bibitem{xie2019scalable}
Y.~Xie, M.~Chen, H.~Jiang, T.~Zhao, and H.~Zha, ``On scalable and efficient computation of large scale optimal transport,'' \emph{arXiv:1905.00158}, 2019.

\bibitem{makkuva2020optimal}
A.~Makkuva, A.~Taghvaei, S.~Oh, and J.~Lee, ``Optimal transport mapping via input convex neural networks,'' in \emph{International Conference on Machine Learning}.\hskip 1em plus 0.5em minus 0.4em\relax PMLR, 2020, pp. 6672--6681.

\bibitem{korotin2021wasserstein}
A.~Korotin, V.~Egiazarian, A.~Asadulaev, A.~Safin, and E.~Burnaev, ``{W}asserstein-2 generative networks,'' in \emph{International Conference on Learning Representations}, 2021.

\bibitem{hua2013implementation}
M.-D. Hua, G.~Ducard, T.~Hamel, R.~Mahony, and K.~Rudin, ``Implementation of a nonlinear attitude estimator for aerial robotic vehicles,'' \emph{IEEE Transactions on Control Systems Technology}, vol.~22, no.~1, pp. 201--213, 2013.

\bibitem{barrau2014intrinsic}
A.~Barrau and S.~Bonnabel, ``Intrinsic filtering on {L}ie groups with applications to attitude estimation,'' \emph{IEEE Transactions on Automatic Control}, vol.~60, no.~2, pp. 436--449, 2014.

\bibitem{barczyk2015invariant}
M.~Barczyk, S.~Bonnabel, J.~Deschaud, and F.~Goulette, ``Invariant {EKF} design for scan matching-aided localization,'' \emph{IEEE Transactions on Control Systems Technology}, vol.~23, no.~6, pp. 2440--2448, 2015.

\bibitem{hesch2014camera}
J.~A. Hesch, D.~G. Kottas, S.~L. Bowman, and S.~I. Roumeliotis, ``Camera-imu-based localization: Observability analysis and consistency improvement,'' \emph{The International Journal of Robotics Research}, vol.~33, no.~1, pp. 182--201, 2014.

\bibitem{kwon2013geometric}
J.~Kwon, H.~S. Lee, F.~C. Park, and K.~M. Lee, ``A geometric particle filter for template-based visual tracking,'' \emph{IEEE transactions on pattern analysis and machine intelligence}, vol.~36, no.~4, pp. 625--643, 2013.

\bibitem{choi2012robust}
C.~Choi and H.~I. Christensen, ``Robust 3d visual tracking using particle filtering on the special euclidean group: A combined approach of keypoint and edge features,'' \emph{The International Journal of Robotics Research}, vol.~31, no.~4, pp. 498--519, 2012.

\bibitem{mccann2001polar}
R.~J. McCann, ``Polar factorization of maps on {R}iemannian manifolds,'' \emph{Geometric \& Functional Analysis GAFA}, vol.~11, no.~3, pp. 589--608, 2001.

\bibitem{ambrosio2005gradient}
L.~Ambrosio, N.~Gigli, and G.~Savar{\'e}, \emph{Gradient flows: in metric spaces and in the space of probability measures}.\hskip 1em plus 0.5em minus 0.4em\relax Springer Science \& Business Media, 2005.

\bibitem{kovachki2020conditional}
R.~Baptista, B.~Hosseini, N.~B. Kovachki, and Y.~Marzouk, ``Conditional sampling with monotone gans: from generative models to likelihood-free inference,'' \emph{arXiv:2006.06755}, 2023.

\bibitem{shi2022conditional}
Y.~Shi, V.~De~Bortoli, G.~Deligiannidis, and A.~Doucet, ``Conditional simulation using diffusion {S}chr{\"o}dinger bridges,'' in \emph{Uncertainty in Artificial Intelligence}.\hskip 1em plus 0.5em minus 0.4em\relax PMLR, 2022, pp. 1792--1802.

\bibitem{spantini2019coupling}
A.~Spantini, R.~Baptista, and Y.~Marzouk, ``Coupling techniques for nonlinear ensemble filtering,'' \emph{SIAM Review}, vol.~64, no.~4, pp. 921--953, 2022.

\bibitem{taghvaei2022optimal}
A.~Taghvaei and B.~Hosseini, ``An optimal transport formulation of {B}ayes’ law for nonlinear filtering algorithms,'' in \emph{2022 IEEE 61st Conference on Decision and Control (CDC)}.\hskip 1em plus 0.5em minus 0.4em\relax IEEE, 2022, pp. 6608--6613.

\bibitem{doucet09}
A.~Doucet and A.~M. Johansen, ``A tutorial on particle filtering and smoothing: Fifteen years later,'' \emph{Handbook of nonlinear filtering}, vol.~12, no. 656-704, p.~3, 2009.

\bibitem{cheng2010particle}
Y.~Cheng and J.~L. Crassidis, ``Particle filtering for attitude estimation using a minimal local-error representation,'' \emph{Journal of Guidance, Control, and Dynamics}, vol.~33, no.~4, pp. 1305--1310, 2010.

\bibitem{vernaza2006rao}
P.~Vernaza and D.~D. Lee, ``Rao-{B}lackwellized particle filtering for 6-{DOF} estimation of attitude and position via {GPS} and inertial sensors,'' in \emph{Proceedings of the IEEE International Conference on Robotics and Automation}, 2006, pp. 1571--1578.

\bibitem{crassidis2003unscented}
J.~L. Crassidis and F.~L. Markley, ``Unscented filtering for spacecraft attitude estimation,'' \emph{Journal of Guidance, Control, and Dynamics}, vol.~26, no.~4, pp. 536--542, 2003.

\bibitem{bonnabel2009IEKF}
S.~Bonnabel, P.~Martin, and E.~Sala{\"u}n, ``Invariant extended {K}alman filter: theory and application to a velocity-aided attitude estimation problem,'' in \emph{Proceedings of the 48th IEEE Conference on Decision and Control}, 2009, pp. 1297--1304.

\bibitem{barrau2015TAC}
A.~Barrau and S.~Bonnabel, ``Intrinsic filtering on {L}ie groups with applications to attitude estimation,'' \emph{IEEE Transactions on Automatic Control}, vol.~60, no.~2, pp. 436--449, 2015.

\bibitem{brenier1991polar}
Y.~Brenier, ``Polar factorization and monotone rearrangement of vector-valued functions,'' \emph{Communications on pure and applied mathematics}, vol.~44, no.~4, pp. 375--417, 1991.

\bibitem{cohen2021riemannian}
S.~Cohen, B.~Amos, and Y.~Lipman, ``Riemannian convex potential maps,'' in \emph{International Conference on Machine Learning}.\hskip 1em plus 0.5em minus 0.4em\relax PMLR, 2021, pp. 2028--2038.

\bibitem{rezende2020normalizing}
D.~J. Rezende, G.~Papamakarios, S.~Racaniere, M.~Albergo, G.~Kanwar, P.~Shanahan, and K.~Cranmer, ``Normalizing flows on tori and spheres,'' in \emph{International Conference on Machine Learning}.\hskip 1em plus 0.5em minus 0.4em\relax PMLR, 2020, pp. 8083--8092.

\bibitem{lou2020neural}
A.~Lou, D.~Lim, I.~Katsman, L.~Huang, Q.~Jiang, S.~N. Lim, and C.~M. De~Sa, ``Neural manifold ordinary differential equations,'' \emph{Advances in Neural Information Processing Systems}, vol.~33, pp. 17\,548--17\,558, 2020.

\bibitem{mathieu2020riemannian}
E.~Mathieu and M.~Nickel, ``Riemannian continuous normalizing flows,'' \emph{Advances in Neural Information Processing Systems}, vol.~33, pp. 2503--2515, 2020.

\bibitem{de2022riemannian}
V.~De~Bortoli, E.~Mathieu, M.~Hutchinson, J.~Thornton, Y.~W. Teh, and A.~Doucet, ``Riemannian score-based generative modeling,'' \emph{arXiv:2202.02763}, 2022.

\bibitem{benamou2000computational}
J.-D. Benamou and Y.~Brenier, ``A computational fluid mechanics solution to the {M}onge-{K}antorovich mass transfer problem,'' \emph{Numerische Mathematik}, vol.~84, no.~3, pp. 375--393, 2000.

\bibitem{evensen2006}
G.~Evensen, \emph{Data Assimilation: The Ensemble {K}alman Filter}.\hskip 1em plus 0.5em minus 0.4em\relax Springer Science \& Business Media, 2006.

\bibitem{Zhu2018}
A.~Z. Zhu, D.~Thakur, T.~\"Ozaslan, B.~Pfrommer, V.~Kumar, and K.~Daniilidis, ``The multivehicle stereo event camera dataset: An event camera dataset for 3d perception,'' \emph{IEEE Robotics and Automation Letters}, vol.~3, no.~3, pp. 2032--2039, 2018.

\end{thebibliography}

\end{document}